\documentclass[twoside,11pt]{article}

\usepackage{blindtext}
\usepackage{amsmath,amsthm}
\usepackage{algorithm}%

\usepackage{graphicx}
\usepackage{epstopdf}
\usepackage{algorithmicx}%
\usepackage{algpseudocode}%
\usepackage{listings,comment}%
\newtheorem{assumption}{Assumption}%
\usepackage{jmlr2e}

\usepackage{lastpage}
\jmlrheading{23}{2022}{1-\pageref{LastPage}}{1/21; Revised 5/22}{9/22}{21-0000}{Author One and Author Two}

\ShortHeadings{Sample JMLR Paper}{One and Two}
\firstpageno{1}

\begin{document}

\title{Delayed supermartingale convergence lemmas for stochastic approximation with Nesterov momentum}

\author{\name Zhang Ming-Kun \email mkzhang@shzu.edu.cn\\
       \addr Colledge of Sciences\\
       Shihezi University\\
       Xinjiang, China}

\editor{My editor}

\maketitle

\begin{abstract}%   <- trailing '%' for backward compatibility of .sty file
 This paper focus on the convergence of stochastic approximation with Nesterov momentum. Nesterov acceleration has proven effective in machine learning for its ability to reduce computational complexity. The issue of delayed information in the acceleration term remains a challenge to achieving the almost sure convergence. Based on the delayed supermatingale convergence lemmas, we give a series  of framework for almost sure convergence. Our framework applies to several widely-used random iterative methods, such as stochastic subgradient methods, the proximal Robbins-Monro method for general stochastic optimization, and the proximal stochastic subgradient method for composite optimization. Through the applications of our framework, these methods with Nesterov acceleration achieve almost sure convergence. And three groups of numerical experiments is to check out theoretical results.
\end{abstract}

\section{Introduction}

Stochastic approximation methods have gained significant prominence in addressing optimization challenges across diverse fields, particularly in the context of machine learning and risk management.  The algorithms as stochastic gradient descent (SGD) \cite{Robbins1951A} and proximal Robbins-Monro methods \cite{Toulis2020The} are well-regarded for their efficiency and memory cost. However, achieving convergence, especially for methods without an inherent delay mechanism, presents a significant challenge. In recent years, stochastic iterative methods have become notable contenders for addressing optimization issues, specially when dealing with large datasets. Moreover, the incorporation of acceleration techniques such as Nesterov's momentum has further enhanced the efficiency and convergence speed of stochastic approximation methods, underscoring their widespread acknowledgment and applicability in various optimization packages, including \verb|keras.optimizers|, \verb|paddle.optimizer|, \verb|sklearn.neural_network|, and \verb|torch.optim|.

\subsection{Stochastic approximation methods}

 The stochastic optimization problem can be formulated as minimizing the expected value of a function,  denoted by $f(x) = \mathbb{E}[F(x, \xi)]$, where $F$ depends on both the decision variable $x$ and a random variable $\xi$. In the context of machine learning, the sheer volume of data samples often renders direct calculation of this expectation computationally expensive. Similarly, risk management, where $\xi$ might represent a continuous random variable, faces analytical intractability when determining the expectation.

 Stochastic approximation methods have emerged as powerful tools to address these challenges. By iteratively sampling from the underlying distribution, these methods estimate the expected value, effectively circumventing the computational hurdles imposed by massive datasets or complex random variable structures. Consequently, stochastic approximation methods have become ubiquitous across diverse fields, including machine learning and risk management. However, the convergence rate of stochastic approximation method is slow,  which  the modified method is required.

\subsection{Stochastic approximation methods with momentum}

To further enhance the efficiency and convergence speed of stochastic approximation methods, techniques like Nesterov acceleration have been incorporated. Nesterov's method, built on the concept of momentum-driven optimization, leverages gradient information more effectively by incorporating a smoothing component that accounts for past updates \cite{assran2020convergence}. This integration has demonstrably improved the performance of stochastic approximation algorithms, primarily by reducing the number of iterations required to achieve a desired level of accuracy. However, it's important to note that the effectiveness of momentum-based methods, like Nesterov acceleration, can be sensitive to the chosen parameters. For instance, using a suboptimal momentum value in Nesterov's method might not outperform the original stochastic gradient descent \cite{Rahul2018The}.

The stochastic approximation algorithm enhanced with Nesterov acceleration operates as follows:
$$
\begin{aligned}
\textbf{(Step 1)}~ & v_{k+1} =x_k - \alpha_k g(x_k,\xi_k)\\
\textbf{(Step 2)}~ & x_{k+1} = (1+\theta_k)v_{k+1} - \theta_k v_k
\end{aligned}
$$
where $g(x_k,\xi_k)$ denotes stochastic first-order information of the objective function,
$\alpha_k>0$ is the step size, and $\theta_k$ is the momentum parameter. When $\theta_k\equiv 0$,  it reverts to the conventional stochastic approximation approach. For $\theta_k\in (-1,0)$,  it signifies a weighted delayed stochastic approximation variant,  while $\theta_k\in (0,1)$ corresponds to the Nesterov accelerated stochastic approximation method,  which is the focus of this paper. $v_k$ could be seen as a delay term for $x_{k+1}$.

Existing convergence analysis, such as the Robbins-Siegmund lemma \cite{Robbins1971A}, have played a crucial role in establishing the almost sure convergence of other stochastic iterative methods. However, its reliance on a specific analytical framework limits its applicability to methods lacking inherent delays. Alternative approaches, such as the dynamical system perspective \cite{Benaim2005Stochastic}, and recent advancements like the composite stochastic optimization coupling supermartingale and T-coupling Supermartingale \cite{Wang2017Stochastic}, \cite{Yang2019Multilevel}, offer promising avenues for further analysis.

Previous works have extensively analyzed the efficiency of Nesterov acceleration under various conditions, including quadratic objectives \cite{assran2020convergence}, \cite{Safavi2018An} and smooth, strongly convex functions \cite{Jain2018Accelerating}. In the presence of smooth and strongly convex conditions, the almost sure convergence rate has been elucidated for methods with a constant momentum parameter \cite{Liu2022On}. Additionally, the convergence of the expected function value has been demonstrated for smooth and nonconvex functions \cite{Liang2023Stochastic}. Further research has introduced dynamical strategies for adjusting the momentum parameter \cite{Sun2021A}, \cite{Sun2022Scheduling}, while other works have proposed quasi-hyperbolic momentum related to Nesterov momentum and analyzed its almost sure convergence under smooth conditions \cite{Zhou2020Amortized}. Finally, the role of memory in stochastic optimization has also been discussed in the literature \cite{Gitman2019Understanding}.

\subsection{Our Contributions}

\begin{itemize}
\item \textbf{Novel supermartingale convergence lemmas with delay:} This paper aims to bridge the gap in the existing literature by introducing a novel approach that reveals a framework for constructing custom supermartingale sequences tailored to specific stochastic iterative methods. Unlike the classical stochastic approximation,  the almost sure convergence of stochastic approximation with Nesterov momentum require different versions of the supermatingale convergence lemmas.

\item \textbf{Almost sure convergence for stochastic subgradient method with Nesterov momentum}  By leveraging the "delay" structure, we provide a novel insight into the almost sure convergence for Nesterov accelerated stochastic approximation without differentiable assumptions. Projection operator onto convex set is also allowed.

\item \textbf{Almost sure convergence for proximal methods with Nesterov momentum} Nesterov accelerated proximal Robbins-Monro methods obtains almost sure convergence. For composite optimization, proximal stochastic subgradient method with Nesterov momentum also obtains almost sure convergence.
\end{itemize}

The rest of this paper is organized as follows. In section 2, a convergence lemma for supermartingales with delay term is presented, which serves as the theoretical foundation for proving the convergence of Nesterov's accelerated stochastic gradient method. In section 3 and section 4, the stochastic subgradient method with Nesterov acceleration and the proximal Robbins-Monro method with Nesterov acceleration obtains almost surely convergence.  In section 5, the stochastic proximal gradient method for composite optimization obtains almost surely convergence. In the last section, we give some numerical experiments for these methods.

\section{A supermartingale convergence lemma with delay}
\label{sec2}

In this section, we explore the convergence guarantees of stochastic processes for the stochastic approximation methods with momentum. We start by examining the stability and convergence of a matrix systems, where the boundedness of the infinite product of matrices $\{M_k\}$ (Proposition \ref{prop:M}) and the convergence of companion matrices associated with quadratic polynomials (Proposition \ref{prop_converge}) provide foundational results.

These matrix-focused propositions are complemented by a fixed-point characterization (Proposition \ref{prop_fix}) that elucidates the limiting behavior of the optimization process, offering a deeper understanding of the algorithm's long-term dynamics. The boundedness and monotonic properties of a related sequence of matrices $\{Q_n\}$ (Proposition \ref{coro}) further reinforce the convergence guarantees of the system. Proposition \ref{prop_stop} presents a critical result ensuring the almost sure convergence of a stochastic sequence $r_n$ to a random variable $V_\infty$, conditional on the convergence of another stochastic sequence $V_n$ and subject to a specific recursive relationship governed by the sequence $\theta_n$. This proposition is instrumental in analyzing the convergence of optimization algorithms under noise.

In this section, the cornerstone of our convergence analysis is Lemma \ref{lemma_Doob}, which offers a vital inequality for bounding the expected future values of a nonnegative stochastic process. This lemma holds particular relevance for the convergence analysis of the Nesterov accelerated method, alongside other iterative optimization techniques that function within stochastic settings. And the Lemma \ref{lem_const}  gives the convergence analysis for nesterov acceleration with constant momentum parameter.

\subsection{The stability of a second-order difference equation}

The critical challenge within this methodological framework stems from Step 2. It is fascinating to observe that the second-order difference equation shares a similar structural design, as exemplified by the equation:
$$ r_{k+1} = a_{k,1} r_k + a_{k,2} r_{k-1}.$$
This construction can be translated into a matrix system devoid of delay,
as represented by the equation:
$$\rho_{k+1} = M_k \rho_k,$$
where $\rho_k = \left[\begin{matrix} r_{k-1} \\ r_{k} \end{matrix}\right]$,
and $M_k = \left[\begin{matrix} 0 & a_{k,1} \\ 1 & a_{k,2} \end{matrix} \right]$.
The stability of such a second-order stochastic difference equation hinges
upon the characteristics of the infinite matrix product $\prod_{k=1}^\infty M_k \triangleq M_1 M_2 M_3\cdots$.
Consequently, in this section, we articulate propositions concerning the infinite production of matrices.
Proposition \ref{prop:M} elucidates the boundedness properties of the infinite sequence of matrices.
Proposition \ref{prop:M} offers a sufficient condition
for the convergence of the infinite matrix sequence derived from Step 2.
Lastly, Proposition \ref{prop_fix} delineates the construction of a fixed-point for the matrix system,
which is inherently a solution to the infinite production.

\begin{proposition}
\label{prop:M}
Let $\{M_k\}$ be a sequence of $2 \times 2$ matrices.
$M_k$ has eigenvalues $\{1, \lambda_k\}$ with $\lambda_k \in (-1, 1)$ for all $k \geq 1$.
Then the infinite product of these matrices, $\prod_{k=1}^\infty M_k$,
converges in the spectral norm, i.e., $\left\| \prod_{k=1}^\infty M_k \right\|_2 < \infty$.
\end{proposition}
\begin{proof}
Let $v$ be an arbitrary vector in $\mathbb{R}^2$. For each $k \geq 1$, it follows from the eigenvalue condition that:
$$
v^\top M_k^\top M_k v \leq v^\top v.
$$
This inequality implies that the spectral norm of $M_k$ is bounded by $1$, i.e., $\|M_k\|_2 \leq 1$.  By mathematical induction, we establish that for any $n \geq 1$:
$$ \left\| \prod_{k=1}^n M_k v \right\|_2 \leq \|v\|_2, $$
where $\prod_{k=1}^n M_k= M_1 M_2 \cdots M_n$.
Since this inequality holds for all $v \in \mathbb{R}^2$, it follows that:
$$\left\| \prod_{k=1}^n M_k \right\|_2 \leq 1,$$
for all $n \geq 1$. As the spectral norm is a continuous function, we can take the limit as $n \rightarrow \infty$ to obtain:
$$ \left\| \prod_{k=1}^\infty M_k \right\|_2 \leq 1 < \infty.$$
\end{proof}

\begin{proposition}
\label{prop_converge}
Consider the companion matrix $M_k$ of the polynomial of degree $2$,
$$
P_k(x)=(x-1)(x-\theta_k), \quad k=1,2,\ldots,
$$
where $\theta_k \in (-1,1)$ for all $k \geq 2$ and $\prod_{k=1}^\infty \theta_k = 0$.
The infinite product of matrix sequence $\prod_{k=1}^\infty M_k$ converges.
\end{proposition}

\begin{proof}
The companion matrix of $P_k(x)$ is given by

$$
M_k = \begin{bmatrix}
0 & -\theta_k \\
1 & 1 + \theta_k \\
\end{bmatrix}.
$$

By mathematical induction, we have

$$
P_{n+1} = \prod_{k=1}^{n+1} M_k = \begin{bmatrix}
-\sum_{k=1}^{n} \prod_{j=1}^k \theta_j & -\sum_{k=1}^{n+1} \prod_{j=1}^{k} \theta_j \\
1 + \sum_{k=1}^{n} \prod_{j=1}^k \theta_j & 1 + \sum_{k=1}^{n+1} \prod_{j=1}^{k} \theta_j \\
\end{bmatrix}.
$$

Hence,

$$
P_{n+1} - P_{n} = \prod_{k=1}^{n} \theta_k \begin{bmatrix}
-1 & -\theta_{n+1} \\
1 & \theta_{n+1} \\
\end{bmatrix},
$$
and
$$
\|P_{n+1} - P_{n}\|_F \leq 2\prod_{j=1}^{n} |\theta_j|.
$$
Therefore, $\{P_{n}\}$ is a Cauchy sequence in the Frobenius norm, as
$$
\|P_{n+m} - P_{n+1}\|_F \leq \sum_{k=1}^m \|P_{n+k+1} - P_{n+k}\|_F \leq \sum_{k=1}^m \prod_{j=1}^{n+k} |\theta_j| = \prod_{j=1}^n |\theta_j| \left(\sum_{k=1}^m \prod_{j=1}^k |\theta_{n+j}|\right).
$$
Then matrix sequence $\{P_n\}$ is convergent.
\end{proof}

Notice that the assumption of $\prod_{k=1}^\infty \theta_k$ allows constant sequence
$\theta_k\equiv c\in (0,1), k\geq 1$ and also asymptotic sequence $\theta_k=\frac{1}{k^s}, (s>0)$.

\begin{proposition}
\label{prop_fix}
$\forall t$, $\bar{X}(t)=\left[\begin{matrix}
-t &-t\\
1+t & 1+t
\end{matrix}
\right]$
is a fixed point of system $X_{k+1} =X_k M_k$,
 which means $\bar{X}(t)=\bar{X}(t)M_k, \forall k\geq2$.
 Under the condition of Proposition \ref{prop_converge},
  $\{X_k\}$ converges to some point as $\bar{X}(t)$.
\end{proposition}
\begin{proof}
By the matrix multiplication, the fixed point is obviously,
$$\left[\begin{matrix}
-t & -t\\
1+t& 1+t
\end{matrix}\right] \left[\begin{matrix}
0 &-\theta_k\\
1& 1+\theta_k
\end{matrix}\right]
=\left[\begin{matrix}
-t & -t\\
1+t& 1+t
\end{matrix}\right].$$
Then we discuss the convergence. Denote set $S= \left\{\left[\begin{matrix}
-t & -t\\
1+t &1+t
\end{matrix}
\right], t\in \mathbb{R}\right\}$.

Under the assumption of Proposition \ref{prop_converge}, $P_n$ could be formulated as $\left[\begin{matrix}
-d_n & -c_n\\
1+d_n & 1+c_n
\end{matrix}\right] $,
$$ \mathrm{dist}^2(P_n,S)= \inf_{X\in S} \|P_n-S\|^2_F= (d_n-c_n)^2,$$
where $S^*_n=\Pi_{S}(P_n)= \left[\begin{matrix}
-\frac{d_n+c_n}{2} & -\frac{d_n+c_n}{2}\\
1+\frac{d_n+c_n}{2}& 1+\frac{d_n+c_n}{2}
\end{matrix}\right]$,
and
$$\mathrm{dist}^2(P_{n+1},S) =\mathrm{dist}^2(M_{n+1}P_n,S)= \theta_{n+1}^2 (d_n-c_n)^2.$$
where $S_{n+1}^* = \left[\begin{matrix}
 -\frac{d_n+c_n}{2}\theta_{n+1}  &  -\frac{d_n+c_n}{2} \theta_{n+1}\\
1+\frac{d_n+c_n}{2}\theta_{n+1} &  1+\frac{d_n+c_n}{2}\theta_{n+1}
\end{matrix}\right]$.
Hence
$$\lim_{n\to \infty}\mathrm{dist}^2(P_{n},S)= \lim_{n\to\infty} \prod_{n=1}^\infty \theta_n^2=0.$$
The proposition is proved.
Furthermore,
$P_\infty=\left[\begin{matrix}
-\sum_{k=1}^\infty \prod_{j=1}^k \theta_j & -\sum_{k=1}^\infty \prod_{j=1}^k\theta_j\\
1+\sum_{k=1}^\infty \prod_{j=1}^k \theta_j & 1+\sum_{k=1}^\infty \prod_{j=1}^k\theta_j
\end{matrix}\right]$.

\end{proof}
Based on the Proposition \ref{prop_fix}, it implies Proposition \ref{coro}.
Proposition \ref{coro} extends the convergence of the infinite product of matrices
from starting with the first index to starting with any index.
And establish the relationship of the monotonicity between the sequence $t_n$ and $\theta_n$.
\begin{proposition}
\label{coro}
Denote $$Q_n=\prod_{k=n}^\infty M_k=\left[\begin{matrix}
-t_n &-t_n\\
1+t_n & 1+t_n
\end{matrix}\right],$$
where $\theta_n\in [c,d]\subset [0,1)$, $t_n=\sum_{k=n}^\infty \prod_{j=n}^k \theta_j\geq 0$ and $Q_n=M_{n}Q_{n+1}, \forall n\geq 1$, which means $t_n=(1+t_{n+1})\theta_{n}.$

Then
\begin{itemize}
\item $\{Q_n\}$ is bounded.
\item If $\{\theta_n\}$ is non-increasing sequence, sequence $\{t_n\}$ is also non-increasing.
\end{itemize}
\end{proposition}

 Until now the stability of the second-order difference equation is ready for the "cornerstone" lemma.

\subsection{A supermartingale convergence Lemma with delay}

In the context of stochastic approximation methods without delay,
the Robbins-Siegmund lemma is a powerful tool that can
be used to establish almost sure convergence for a variety of algorithms,
including the stochastic subgradient method, the proximal Robbins-Monro method,
and even some derivative-free methods.
It is observed that the almost sure convergence
for delayed stochastic approximation methods cannot be trivially extended from the non-delayed case.

The presence of delay introduces additional complexities
that necessitate different approaches to ensure almost sure convergence.
Therefore, in this section,
we will explore various versions of the Nesterov accelerated strategy and the weighted average strategy,
which are designed to handle delayed stochastic approximation methods more effectively.

\begin{proposition}
\label{prop_stop}
If stochastic sequence $V_n$ converges to a random variable $V_\infty$ almost surely.
And stochastic sequence $r_n$ satisfies
$$r_{n+1}= (1-\theta_n) r_n + \theta_n V_{n+1}, a.s, \theta_k \in [c,d]\subset [0,1).$$
Then  $r_n$ converges to $V_\infty$, a.s.
\end{proposition}
\begin{proof}
Stochastic sequence $V_n$ converges to a random variable $V_\infty$ almost surely. Set a Markov time $\tau_1$ satisfying
$$\mathcal{P}\left(\omega, \forall \varepsilon>0, \exists \tau_1(\omega)\in (0,+\infty)
, \forall n>\tau_1(\omega), |V_{n}(\omega)-V_\infty(\omega)|<\frac{\varepsilon}{2}\right)=1,$$
where $\tau_1$ is adapted to $\mathcal{F}_n\supset \sigma(V_1,\cdots,V_n)$. Then $\tau_1$ is a stopping time.

According to $\theta_k\in [c,d]\subset [0,1)$, $\prod_{k=1}^\infty \theta_k\leq \lim_{n\to \infty}d^n=0$. $\forall \varepsilon>0$, $\exists N_2$, $\forall n> N_2$, $\prod_{k=\tau_1}^n \theta_k\leq d^{n-\tau_1}<\frac{\varepsilon}{2}$, a.s.  Then $\tau_2=\max\left\{\tau_1, \log_d\left(\frac{\varepsilon}{2\tau_1}\right)+\tau_1\right\}=\tau_1+\log_d\left(\frac{\varepsilon}{2\tau_1}\right)$ is also a stopping time, with $\forall n> \tau_2$,
$$|r_n-V_\infty|<\frac{\varepsilon}{2}+\frac{\varepsilon}{2}=\varepsilon, a.s.$$
\end{proof}

Proposition \ref{prop_stop} gives the relationship of the convergence of $\{r_n\}$ and $\{V_n\}$,
which is important for Lemma \ref{lemma_Doob}.

\begin{lemma}
\label{lemma_Doob}
If stochastic process
$\{r_n\}$ is nonnegative and sequence $\theta_n\in [c,d]\subset [0,1)$, with $\sup_k |r_k|<+\infty$.

$$\mathbb{E}[r_{n+2}|\mathcal{F}_{n+1}]\leq (1+\theta_{n}) r_{n+1} - \theta_{n}r_n,~\theta_n\in [c,d]\subset [0,1),n\geq 1$$

Then $r_n$ converges to some finite random variable $r_\infty$,a.s.
\end{lemma}
\begin{proof}
Set $V_n = \rho_n^\top Q_n \phi$, where $Q_n=\prod_{k=n}^\infty M_k$ in Proposition \ref{coro},
$\rho_n = \begin{pmatrix} r_{n} \\ r_{n+1} \end{pmatrix}$,
$\phi = \begin{pmatrix} \phi_1 \\ \phi_2 \end{pmatrix}$,  $\phi_1+\phi_2 > 0$.
According to Corollary \ref{coro}, $\|Q_n\|_2$ is bounded,
which implies $\{V_n\}$ is bounded.  Furthermore, we have
$$\begin{aligned}
&V_{n+1}\\
=& \left[r_{n+1}, r_{n+2}\right]\left[
\begin{matrix}
-t_{n+1} & -t_{n+1}\\
1+t_{n+1} & 1+t_{n+1}
\end{matrix}\right]\left[\begin{matrix}\phi_1\\ \phi_2\end{matrix}\right]\\
=&(\phi_1+\phi_2)((1+t_{n+1})r_{n+2}-t_{n+1}r_{n+1}).
\end{aligned}.$$

$$\mathbb{E}[V_{n+1}|\mathcal{F}_n]\leq (\phi_1+\phi_2)((1+t_{n+1})((1+y_n)r_{n+1}-y_nr_n)- t_{n+1}r_{n+1})=V_n,$$
This indicates that $t_{n+1}=\sum_{k=n+1}^\infty \prod_{j=k}^\infty\lambda_j>0$.
Hence $V_n$ is a supermartingale.  Since $\{V_n\}$ s bounded and a supermartingale,
it converges to some random variable $V_\infty$ almost surely by Doob' s martingale convergence theorem.
According to Proposition \ref{prop_stop} and $r_{n+1}=\frac{1}{1+y_n}V_{n+1}+\frac{y_n}{1+y_n}r_n, n\geq 1$, a.s, $r_n$ converges to $V_\infty$, a.s. $\sum_{k=1}^\infty \eta_k<+\infty$.
\end{proof}

Lemma \ref{lemma_Doob} is a basic version as the Doob's Submartingale convergence theorem for stochastic approximation method without delay, which is essential to all the following lemmas. Yet, it remains challenging to establish a supermartingale with the required lower boundedness for supermartingale with lower boundedness in Nesterov accelerated methods. Fortunately, under the assumption of uniform boundedness of the iterative points, we can also give the almost surely convergence.

If the momentum $\theta_k\equiv \theta\in (0,1)$ is a constant, $t_n\equiv t= \frac{\theta}{1-\theta}$.
If $(1+\theta) r_{n+1}-\theta r_n\geq 0$,  $(1+t)r_{n+1}-t r_n= \frac{1}{1-\theta}((1+\theta)r_n-\theta r_n)\geq 0$.
Then $V_n$ is nonnegative supermartingale, which converges almost surely.
Hence the boundedness of $r_k$ can be removed.

\begin{lemma}
\label{lem_const}
If a stochastic process $\{r_n\}$ is nonnegative and sequence $\theta_n\equiv \theta\in (0,1)$.
$$\mathbb{E}[r_{n+2}|\mathcal{F}_{n+1}]\leq (1+\theta) r_{n+1} - \theta r_n,n\geq 1.$$
Then $r_n$ converges to some finite random variable $r_\infty$, a.s.
\end{lemma}

\section{Application in Stochastic subgradient methods with Nesterov acceleration}

In this section, we mainly consider the stochastic subgradient methods (ssgd) with Nesterov acceleration for the simple set constrained stochastic optimization,

$$\min_{x\in C} f(x)=\mathbb{E}[F(x,\xi)],$$
where $C\subset \mathbb{R}^n$ is a convex set.

\begin{algorithm}[H]
\caption{The ssgd method} \label{algo_grad}
\begin{algorithmic}[1]
\Require Step size $\{\alpha_k\}$, momentum size $\{\theta_k\}$, initial value $v_1$, $v_2$,
\For {$n=1,2,\cdots$}
Calculate the Nesterov acceleration,
$$x_{k+1}=(1+\theta_k)v_k-\theta_kv_{k-1}.$$
 Generate a random variable $\xi_{k+1}$ and calculate a subgradient $g(x_{k+1},\xi_{k+1})\in \partial_x F(x_{k+1},\xi_{k+1})$.
$$v_{k+1}= \Pi_C(x_{k+1}-\alpha_k g(x_{k+1},\xi_{k+1})).$$
\EndFor
\end{algorithmic}
\end{algorithm}
$\Pi_C(\cdot)$ is the projection operator to convex set $C$. When $C=\mathbb{R}^n$, the algorithm is known as the popular NAG-SGD method.

\subsection{Lemmas for delayed stochastic subgradient methods}

Lemma \ref{lemma_RS} build a supermartingale with second-order delayed random variable. It gives a common frame for delayed SA methods.  By the arbitrary of positive $\phi_1,\phi_2$, we set $\phi_1+\phi_2=1$.

\begin{lemma}
\label{lemma_RS}
If stochastic process $\{r_n\}$ is nonnegative and sequences  $\{\beta_n\}$, $\{\eta_n\}$ are positive, with $\sup_{n\geq 1} |r_n|<+\infty$.
$$\mathbb{E}[r_{n+2}|\mathcal{F}_{n+1}]\leq \left((1+y_{n}) r_{n+1} - y_{n}r_n\right)+\beta_n-\eta_n, n\geq 1$$
\begin{equation}
\label{pos_asm}
y_n\in [c,d]\subset (0,1),n\geq 1,
\end{equation}
$$\sum_{k=1}^\infty \beta_n<\infty,a.e.$$
Then $r_n$ converges to some finite random variable almost surely, $\sum_{k=1}^\infty \eta_k<+\infty$ almost surely.
\end{lemma}
\begin{proof}
Consider the stochastic sequence defined as:

$$V_n = \rho_n^\top Q_n \phi + 2\sum_{k=n}^\infty \beta_k.$$

 According to $\sum_{k=1}^\infty \alpha_k < \infty$, the sequence $\|Q_n\|_2< \infty$ almost surely.   Finally, with $\sum_{n=1}^\infty \beta_n < \infty$, we conclude that $V_n$ is bounded almost surely, $\forall n\geq 1$

Take the condition expectation on the $\sigma$-algebra $\mathcal{F}_n$,

$$
\begin{aligned}
\mathbb{E}[V_{n+1}|\mathcal{F}_{n+1}]&=
[\mathbb{E}[\rho_{n+1}^\top|\mathcal{F}_{n+1}] ] Q_{n+1}\phi+(\phi_1+\phi_2)\sum_{k=n+1}^\infty \beta_k\\
&\leq [r_{n+1},\quad ((1+y_n) r_{n+1}-y_nr_n)+\beta_n]   Q_{n+1} \phi\\
&\quad +2 (\phi_1+\phi_2)\sum_{k=n+1}^\infty \beta_k\\
&\leq  [r_{n+1}+\beta_n,\quad (1+y_n)r_{n+1}-y_nr_n+\beta_n  ]\left[ \begin{matrix}
-t_{n+1}\\1+t_{n+1}
\end{matrix}\right] \\
&\quad +2\sum_{k=n+1}^\infty \beta_k + t_{n+1} \beta_n\\
&\leq\rho_n^\top Q_{n} \phi + 2\sum_{k=n}^\infty \beta_k\\
&=V_n.\\
\end{aligned}$$

$V_n$ is a nonnegative supermartingale and converges to some finite variable $V_\infty$ almost surely. And $\mathbb{E}[{r}_{n+2}|\mathcal{F}_{n+1}]=(1+t_{n}) r_{n+1}-t_n r_n$, $r_{n+1}\geq \frac{1}{t_n+1}\mathbb{E}[{r}_{n+2}|\mathcal{F}_{n+1}]+\frac{t_n}{1+t_n+1} r_n\geq 0$. Then $V_n$ converges almost surely.
And according to Lemma \ref{lemma_Doob}, $r_n$ converges to $V_\infty$ almost surely.
\end{proof}

\begin{lemma}
If stochastic process $\{r_n\}$ and $\{z_n\}$ are nonnegative and sequences $\{\beta_n\}$, $\{\eta_n\}$ are positive.
\label{lemma_ex}
$$\mathbb{E}[r_{n+2}+z_{n+2}|\mathcal{F}_{n+1}]\leq \left((1+\theta_{n}) r_{n+1} - \theta_{n}r_n+ z_{n+1}\right)+\beta_n-\eta_n$$
\begin{equation}
y_n\in (0,1),  n\geq 2
\end{equation}
is a non-increasing sequence.
$$\sum_{k=1}^\infty \beta_n<\infty,a.e.$$
Set $\rho_n=\left[\begin{matrix} r_n\\ r_{n+1}+ z_{n+1}\end{matrix}\right]$, $Q_n=\left[
\begin{matrix}
-t_{n} & -t_{n}\\
1+t_{n} & 1+t_{n}
\end{matrix}\right]$, $\phi=\left[
\begin{matrix}
\phi_1\\
\phi_2
\end{matrix}\right]$. $$V_n = \rho_n^\top  Q_n \phi +  2\sum_{k=n}^\infty \beta_k.$$
Then $V_n$ converges to some finite random variable almost surely. $\sum_{k=1}^\infty \eta_n<\infty$  almost surely, $z_{n}$ converges to 0 almost surely and $r_{n}$ converges to some random variable $r_\infty$ almost surely.
\end{lemma}
\begin{proof}
By the Corollary \ref{coro},
$$\mathbb{E}[V_{n+1}|\mathcal{F}_n]-V_n\leq (t_{n+1}-t_n)z_n\leq 0.$$
$V_n$ is a bounded supermartingale, which converges to some random variable almost surely.
%By  $\lim\limits_{n\to \infty} \theta_n=0$, $t_n = \theta_n(1+t_{n+1})\to 0$, $n\to \infty$, where $t_n$ is bounded by Proposition \ref{prop:M} and Corollary \ref{coro}. Then
%$$V_\infty=\lim\limits_{n\to \infty} (1+t_n) r_{n+1}-t_n r_n+z_n=\lim_{n\to \infty} r_{n+1}+z_{n+1}.$$
\end{proof}

The following is a kind of coupling version of supermatingale convergence result of Lemma \ref{lemma_ex}.

\begin{lemma}
\label{lemma_couple}
Consider two stochastic sequence $\{r_n\}$, $\theta_n\in [c,d]\subset [0,1)$.
\begin{equation}
\begin{aligned}
&\mathbb{E}[r_{n+2}|\mathcal{F}_{n+1}] \leq ((1+\theta_n) r_{n+1}-\theta_n r_n) -\eta_n + \beta_n + h \zeta_n z_{n+1}\\
&\mathbb{E}[z_{n+2}|\mathcal{F}_{n+1}] \leq (1-\zeta_n) z_{n+1} - \bar{\eta}_{n} + \bar{\beta}_n .
\end{aligned}
\end{equation}
Then $\{r_n\}$ and $\{z_n\}$ converge to some finite variable, a.s and $\sum_{k=1}^\infty \eta_k <+\infty$, a.s. Furthermore, if $\sum_{n=1}^\infty \zeta_n=\infty$, $z_{n+1}$ converges to 0 almost surely and $r_{n+1}$ converges to some random variable $r_\infty$ almost surely.
\end{lemma}
\begin{proof}
Take $J_n= r_{n+1} + hz_{n+1}$.
 $$\mathbb{E}[J_{n+1}|\mathcal{F}_{n+1}] \leq  J_n+\theta_n(r_{n+1}-r_n)   - (\eta_n+h\bar{\eta}_n) + (\beta_n+h\bar{\beta}_n)   $$
Then $V_n$ in Lemma \ref{lemma_ex} converges to some finite random variable almost surely.  So $z_n$ are almost surely bounded.  $\sum_{n=1}^\infty \zeta_n z_{n+1}<+\infty$.  Again by Lemma \ref{lemma_ex}, $\sum_{k=1}^\infty \eta_n< +\infty$, a.s.  $\sum_{k=1}^\infty \bar{\eta}_n< +\infty$, a.e. and $\sum_{k=1}^\infty \zeta_nz_n<+\infty$, a.s.
 So $\{r_n\}$ and $\{z_n\}$ converge to some finite variable, a.s.
\end{proof}

\begin{lemma}
    \label{lemma_const_cvg}
Consider two stochastic sequence $\{r_n\}$, $\theta_n\equiv \theta\in (0,1)$.
\begin{equation}
\begin{aligned}
    &\mathbb{E}[r_{n+2}|\mathcal{F}_{n+1}] \leq ((1+\theta_n) r_{n+1}-\theta_n r_n) -\eta_n + \beta_n + h \zeta_n z_{n+1}\\
    &\mathbb{E}[z_{n+2}|\mathcal{F}_{n+1}] \leq (1-\zeta_n) z_{n+1} - \bar{\eta}_{n} + \bar{\beta}_n .
\end{aligned}
\end{equation}
Then $\{r_n\}$ and $\{z_n\}$ converge to some finite variable, a.s and $\sum_{k=1}^\infty \eta_k <+\infty$, a.s. Furthermore, if $\sum_{n=1}^\infty \zeta_n=\infty$, $z_{n+1}$ converges to 0 almost surely and $r_{n+1}$ converges to some random variable $r_\infty$ almost surely.
\end{lemma}

\subsection{Almost surely convergence}

Consider the Nesterov accelerated projected stochastic subgradient method. The first step is the Nesterov acceleration. And the second step is the projected stochastic subgradient method.  When $C$ is a bounded set, sequence $\{x_k\}$ is naturally bounded.

\begin{assumption}
\label{asm_1}
(a) $F(\cdot,\xi)$ is continues convex for almost sure $\xi\in \Xi$. \\
(b) Subgradient $g(\cdot,\xi)$ of $F(\cdot,\xi)$ a.s. $\xi\in \Xi$.\\
(c) Step size $\alpha_k\geq 0$, satisfies $\sum_{k=1}^\infty \alpha_k=\infty$, $\sum_{k=1}^\infty \alpha_k^2<\infty$. Momentum size $\theta_k \in [c,d]\subset (0,1)$.  \\
(d) Iteration $v_k$ is bounded by $M$, $\sup_{k\geq 1} \|x_k\|\leq M$.
\end{assumption}

Futhermore, if the momumtum parameter $\theta$ is a constant, the boundedness of the iteration $\{v_n\}$ could be relaxed to the boundedness of subdifferential at $v_n$, for example Lipschitz continuous function.
\begin{assumption}
    \label{asm_2}
    (a) $F(\cdot,\xi)$ is continues convex for almost sure $\xi\in \Xi$. \\
    (b) Subgradient $g(\cdot,\xi)$ of $F(\cdot,\xi)$ a.s. $\xi\in \Xi$.\\
    (c) Step size $\alpha_k\geq 0$, satisfies $\sum_{k=1}^\infty \alpha_k=\infty$, $\sum_{k=1}^\infty \alpha_k^2<\infty$. Momentum size $\theta_k \in [c,d]\subset (0,1)$.  \\
    (d) The norm of $\partial F(x,\xi)$ is almost surely bounded,
    $$\sup_{x\in X} \|g(x,\xi)\| \leq M, \forall g(x,\xi)\in \partial F(x,\xi),a.e.\xi\in \Xi.$$
\end{assumption}

\begin{theorem}
\label{thm:ssg}
Consider the sequence $x_k$, $v_k$ generated by Algorithm \ref{algo_grad}, with Assumption \ref{asm_1}.  Then the sequences of $\{x_k\}$ and $\{v_k\}$ both converges to the same optimal.
\end{theorem}
\begin{proof}

$$\|x_{k+1}-x^*\|^2 = (1+\theta_k)\|v_k-x^*\|^2 -\theta_k \|v_{k-1}-x^*\|^2 +\theta_k(1+\theta_k)\|v_k-v_{k-1}\|^2,$$
where by the nonexpansive of projection operator $\Pi_C(\cdot)$
$$\begin{aligned}
\|v_{k}-v_{k-1}\|^2&=\|\Pi_C(x_{k}-\alpha_{k-1} g(x_{k},\xi_{k}))-v_{k-1}\|^2 \\
&\leq \|x_{k}-\alpha_{k-1} g(x_{k,\xi_{k}})-v_{k-1}\|^2\\
&= \|x_k-v_{k-1}\|^2 -2\alpha_{k-1}\langle x_k-v_{k-1},g(x_{k},\xi_k) \rangle+ \alpha_{k-1}^2\|g(x_k,\xi_k)\|^2 \\
\end{aligned}$$

Set $r_{k+1}=\|v_{k+1}-x^*\|^2$ and $z_k=\|v_k-v_{k-1}\|^2$.
Then   $x_k-v_{k-1}$ is replaced by $\theta_{k-1}(v_{k-1}-v_{k-2})$, according to Nesterov acceleration step.
$$\begin{aligned}
 z_k\leq \theta_{k-1}^2\|v_{k-1}-v_{k-2}\|^2 -2\theta_{k-1}\alpha_{k-1}\langle v_{k-1}-v_{k-2},g(x_{k},\xi_k) \rangle+ \alpha_{k-1}^2\|g(x_k,\xi_k)\|^2
\end{aligned}$$

Then by Cauchy-Schwartz inequality $2\langle a,b\rangle\leq \tau \|a\|^2+\frac{1}{\tau}\|b\|^2$, $\forall \tau>0$.
\textsc{$$\begin{aligned}
&\quad z_k\\
& \leq \left(\theta_{k-1}^2+ \frac{\alpha_{k-1}}{\tau_{k-1}} \theta_{k-1} \right)\|v_{k-1}-v_{k-2}\|^2 +(\alpha_{k-1}\theta_{k-1}\tau_{k-1} +\alpha_{k-1}^2)\|g(x_{k},\xi_{k})\|^2\\
&\leq \left(\theta_{k-1}^2+ \frac{\alpha_{k-1}}{\tau_{k-1}} \theta_{k-1} \right)z_{k-1} +(\alpha_{k-1}\theta_{k-1}\tau_{k-1} + \alpha_{k-1}^2)M.
\end{aligned}$$}
Then $z_k$ converges to some finite random variable almost surely by The Robbins-Siegmund Lemma. Without loss of generality, take $\tau_k=\frac{1}{\tau} \alpha_k$.

$$\begin{aligned}
 z_k\leq \left(\theta_{k-1}^2+ \tau \theta_{k-1} \right)z_{k-1} +(\alpha_{k-1}\theta_{k-1}\tau_{k-1} + \alpha_{k-1}^2)M.
\end{aligned}$$

Take $\tau\in (0, \frac{1-d^2}{d})$. Then $p_k=1-\theta_k(\theta_k+\tau)\in (0,1-c^2)$. Take $h\in (0, \frac{1-c^2}{d^2+d})$. $h(\theta_k+\theta_k^2) - p_k \leq 0$.

$$\begin{aligned}
\mathbb{E}[r_{k+1}|\mathcal{F}_k]& \leq \|x_{k+1}-x^*\|^2 -2\alpha_k(f(x_{k+1})-f(x^*))+\alpha_k^2 \mathbb{E}[\|g(x_{k+1},\xi_{k+1})\|^2|\mathcal{F}_{k}]\\
&\leq \|x_{k+1}-x^*\|^2-2\alpha_k(f(x_{k+1})-f(x^*))+O(\alpha_k^2)\\
&\leq
((1+\theta_k)r_{k+1}-\theta_kr_k)-2\alpha_k(f(x_{k+1})-f(x^*))+O(\alpha_k^2)+hp_k z_k
\end{aligned}$$

$$ \mathbb{E}[z_{k+1}|\mathcal{F}_k]\leq (1-p_k)z_k +O(\alpha_k^2). $$

According to $\sum_{k=1}^\infty \alpha^2_k<+\infty$,  $\sum_{k=1}^\infty \theta_k\alpha_k<\infty$.  According to Lemma \ref{lemma_couple} immediately,  $r_k+hz_k$ converges to some finite random variable ,a.s.  $\sum_{k=1}^\infty \alpha_k(f(x_k)-f^*) <\infty$, a.s. There is a subsequence of $x_{k+1}$ converges to $x^*$, a.e. And $\sum_{k=1}^\infty \theta_k(1+\theta_k)\|v_{k+1}-v_k\|^2<\infty$.  By the convexity of function $f$,
$$f(v_{k+1})\leq \frac{1}{1+\theta_k}f(x_{k+1}) +\frac{\theta_k}{1+\theta_k} f(v_k),$$
equally,
$$ - \left( f(x_{k+1})-f^*\right)\leq -\left(f(v_{k+1})-f^*\right)+\theta_k \left(f(v_k)-f(v_{k+1}) \right).$$
Which means there exists a subsequence of $\{\theta_k\}$ converges to 0 according to assumption $\sum_{k=1}^\infty \alpha_k=\infty$. Take $h=1$,
According to $r_k+z_k$ converges almost surely, $\{v_k\}$ and $\{f(v_k)\}$ is bounded,  so

$$
\begin{aligned}
\mathbb{E}[r_{k+1}| \mathcal{F}_k]& \leq (1-2\alpha_k\mu)\left((1+\theta_k)r_k-\theta_{k}r_{k-1}+ \theta_k(1-\theta_k) \|v_k-v_{k-1}\|^2 \right) \\
&\quad +\alpha_k^2M -2\alpha_k(f(v_{k+1})-f(x^*))+  \alpha_k \theta_k M.
\end{aligned}$$
So there exists a subsequence of $\{v_k\}$ converges to $x^*$ and by the convergence of $\|v_k-x^*\|^2$. $v_k$ converges to $x^*$ almost surely.
\end{proof}

\begin{remark}
 The famous parameter sequence $\theta_k= \frac{1}{k+3}$,
 which means $d=\frac{1}{4}$, $\tau \in (0, \frac{15}{4})$,
  $\sum_{k=1}^\infty (1-\theta_k(\theta_k+\tau))=\infty$,
  then the almost surely convergence is obtained.
\end{remark}

According to Lemma \ref{lemma_const_cvg}, Theorem \ref{thm:ssg_cst} is obvious.
\begin{theorem}
    \label{thm:ssg_cst}
Consider the sequence $x_k$, $v_k$ generated by Algorithm \ref{algo_grad} and Assumption \ref{asm_2} holds.  Then the sequences of $\{x_k\}$ and $\{v_k\}$ both converge to the some optimal.
\end{theorem}

\section{Applications in Delayed Proximal Robbins-Monro methods}

Lemma \ref{Lemma_prox0}, \ref{Lemma_prox1} and Lemma \ref{lemma_prox} is prepared for proximal Robbins-Monro method (prox-RM).

With observable noise, proximal Robbin-Monro method with NAG could be displayed as follows.

\textsc{\begin{algorithm}[!htbp]
\caption{The prox-RM method} \label{algo1}
\begin{algorithmic}[1]{
\Require Step size $\{\alpha_k\}$, momentum size $\{\theta_k\}$, initial value $v_1$, $v_2$, and iteration number $N$,
\For {$n=1,2,\cdots$}
Calculate
$$x_{k+1}=(1+\theta_k)v_k-\theta_kv_{k-1},$$
and the proximal point
$$v_{k+1}\in \arg\min_{v\in C} F(v,\xi_k) +\frac{1}{2 \alpha_k} \|v-x_{k+1}\|^2.$$
\EndFor
}
\end{algorithmic}
\end{algorithm}}

\subsection{Lemmas for delayed proximal Robbins-Monro method}
\begin{lemma}
\label{Lemma_prox0}
$\{r_k\}$, $\{\eta_k\}$ are nonnegative stochastic sequences.
$\mathbb{E}[r_{k+1}|\mathcal{F}_k]\leq r_k-a_k (\eta_{k+1}-\eta_k)$,
$a_k\geq 0$ is a decreasing sequence,
Then $r_n+a_{n-1}\eta_n$ converges almost surely to some finite random variable.
\end{lemma}
\begin{proof}
Set $V_n= r_n+a_{n-1}\eta_n\geq 0,a.s.$

$\mathbb{E}[V_{n+1}|\mathcal{F}_k] = r_{n+1}+a_n\eta_{n+1} \leq r_{n} -a_n(\eta_{n+1}-\eta_n)  +a_n \eta_{n+1} \leq r_n +a_{n-1}\eta_n =V_n$.
Then $V_{n}$ converges almost surely.

\end{proof}

\begin{lemma}
\label{Lemma_prox1}
$\{r_k\}$, $\{\eta_k\}$, $\{\beta_k\}$, $\{\zeta_k\}$ are nonnegative stochastic sequences. $\theta_k\in [c,d]\subset [0,1)$ is a decreasing sequence.
$$\mathbb{E}[r_{k+2}|\mathcal{F}_{k+1}]\leq (1+\theta_{k})r_{k+1}-\theta_kr_k-a_k(\eta_{k+1}-\eta_k) +\beta_k-\zeta_k,$$
$a_k\geq 0$ monotone decreasing, $\sum_{k=1}^\infty \beta_k<+\infty$,~a.s.
Then $r_n+a_{n-1}\eta_n$ converges almost surely to some finite random variable.
\end{lemma}
\begin{proof}
Set $V_n= [r_n, r_{n+1}+a_{n}\eta_{n+1}] Q_n\left[\begin{matrix}\phi_1\\ \phi_2\end{matrix}\right]+\sum_{k=n}^\infty \beta_k\geq 0,a.s.$

$$\begin{aligned}
&\quad\mathbb{E}[V_{n+1}|\mathcal{F}_n] \\
&=\mathbb{E}[[r_{n+1}, r_{n+2}+a_{n+1}\eta_{n+2}] |\mathcal{F}_n]Q_n\phi+\sum_{k=n+1}^\infty \beta_k \\
&=\mathbb{E}[[r_{n+1}, (1+\theta_n)r_{n+1}-\theta_nr_n+a_{n}\eta_{n+1}] |\mathcal{F}_n]Q_n\phi+\sum_{k=n+1}^\infty \beta_k \\
&\leq V_n,
\end{aligned}$$
 which is nonnegative supermatingale.
Then $V_{n}$  converges almost surely to some finite random variable. $\sum_{k=1}^\infty \zeta_k<+\infty$.
\end{proof}

\begin{lemma}
\label{lemma_prox}
Consider the  positive stochastic sequence $r_k$, $\eta_k$, $\zeta_k$, $\rho_k$. And sequence $\theta_k\in (0,1)$ is bounded. $a_k$ is a positive decreasing sequence.

$$\mathbb{E}[r_{k+2}|\mathcal{F}_k] \leq (1+\theta_k) r_{k+1}-\theta_k r_k- \eta_k +\beta_k+ hp_k z_{k+1}, k\geq 1$$

$$\mathbb{E}[z_{k+2}|\mathcal{F}_k] \leq (1-p_k) z_{k+1} - a_k(\rho_k-\rho_{k-1}), k\geq 1.$$

Then $r_k$ converges to some finite random variable $r_\infty$,a.s and  $\sum_{k=1}^\infty \eta_k<+\infty$,a.s.
\end{lemma}
\begin{proof}
Set $J_k=r_{k+1}+hz_{k+1}$, the proof is similar to Lemma \ref{lemma_couple}, with Lemma \ref{Lemma_prox1}.
\end{proof}

\begin{lemma}
    \label{lemma_prox_const}
    Consider the  positive stochastic sequence $r_k$, $\eta_k$, $\zeta_k$, $\rho_k$. The momentum parameter $\theta_k\equiv \theta\in  (0,1)$ and the step size $a_k$ is a positive decreasing sequence.

    $$\mathbb{E}[r_{k+2}|\mathcal{F}_k] \leq (1+\theta_k) r_{k+1}-\theta_k r_k- \eta_k +\beta_k+ hp_k z_{k+1}, k\geq 1$$

    $$\mathbb{E}[z_{k+2}|\mathcal{F}_k] \leq (1-p_k) z_{k+1} - a_k(\rho_k-\rho_{k-1}), k\geq 1.$$

    Then $r_k$ converges to some finite random variable $r_\infty$,a.s and  $\sum_{k=1}^\infty \eta_k<+\infty$, a.s.
\end{lemma}

\subsection{Almost surely convergence}

\begin{theorem}
    \label{thm:PRM}
Consider $\{x_k\}$, $\{v_k\}$ generated from Algorithm \ref{algo1} under Assumption \ref{asm_1} and momentum parameter $\{\theta_k\}$ is nonincreasing. Then $\{x_k\}$ and $\{v_k\}$ converges to some optimal $x^*$.
\end{theorem}

\begin{proof}
Assume an arbitrary optimal point $x^*$,
$$\|x_{k+1}-x^*\|^2 = (1+\theta_k) \|v_k-x^*\|^2 - \theta_k \|v_{k-1}-x^*\|^2+ \theta_k(1+\theta_k) \|v_k-v_{k-1}\|^2.$$

$$\|v_{k+1}-x^*\|^2 \leq \|x_{k+1}-x^*\|^2-2\alpha_k (F(v_{k+1},\xi_{k})-F(x^*,\xi_k))- \|v_{k+1}-x_{k+1}\|^2.$$

Both side take conditional expectation on $\sigma$-algebra $\mathcal{F}_{k}$,

$$\begin{aligned}
&\mathbb{E}[\|v_{k+1}-x^*\|^2|\mathcal{F}_k]\\
\leq &\|x_{k+1}-x^*\|^2-2\alpha_k(f(v_{k+1})-f(x^*))- \|v_{k+1}-x_{k+1}\|^2.\\
= & (1+\theta_k) \|v_k-x^*\|^2 - \theta_k \|v_{k-1}-x^*\|^2+ \theta_k(1+\theta_k) \|v_k-v_{k-1}\|^2\\
&-2\alpha_k(f(v_{k+1})-f(x^*))- \|v_{k+1}-x_{k+1}\|^2.
\end{aligned}$$
And
$$\begin{aligned}
&     \qquad \mathbb{E}[\|v_{k+1}-v_k\|^2|\mathcal{F}_k]\\
&\leq \theta_k^2 \|v_k-v_{k-1}\|^2 - \alpha_k(f(v_{k+1})-f(v_k)) - \|v_{k+1}-x_{k+1}\|^2
\end{aligned}.$$
Set $r_k=\|v_k-x^*\|^2$, and $z_k=\|v_k-v_{k-1}\|^2$, $\rho_k=f(v_{k-1})-f^*$.  According to Lemma \ref{lemma_prox} $v_{k}$ converges to some optimal $x^*$, a.s and $\|v_{k+1}-x_{k+1}\|^2\to 0$, a.s. $\|v_k-v_{k-1}\|^2$ converges to 0. $\|v_k-x^*\|$ converges, then converges to 0, a.s.
\end{proof}

\begin{remark}
If the momentum parameter $\theta_k$ is nonincreasing,
Lemma \ref{lemma_prox} could also imply the  convergence of stochastic subgradient method with Nesterov acceleration.
However according to Theorem \ref{thm:ssg}, there is no need for nonincreasing of parameter $\{\theta_k\}$.
\end{remark}

At the end of this section,  the constant momentum parameter case is given.
\begin{assumption}
    \label{asm_3}
    (a) $F(\cdot,\xi)$ is continues convex for almost sure $\xi\in \Xi$. \\
    (b) Subgradient $g(\cdot,\xi)$ of $F(\cdot,\xi)$ a.s. $\xi\in \Xi$.\\
    (c) Step size $\alpha_k\geq 0$, satisfies $\sum_{k=1}^\infty \alpha_k=\infty$, $\sum_{k=1}^\infty \alpha_k^2<\infty$. Momentum size $\theta_k \equiv \theta\in (0,1)$.  \\
\end{assumption}

According to the Lemma \ref{lemma_prox_const}, under the constant momentum parameter, the proximal Robbins-Monro method converges almost surely.

\begin{theorem}
Suppose that $\{v_k\}$, $\{x_k\}$ is generated by the Algorithm \ref{algo2}, and Assumption \ref{asm_3} holds. The $\{x_k\}$ converges to some optimal almost surely.
\end{theorem}

In Assumption \ref{asm_3}, the boundedness of the subgradient and the boundedness of the iterative sequence is removed from Assumption \ref{asm_1} and Assumption \ref{asm_2}.

\section{Application in composite optimization}

Consider the composite optimization as follows,
\begin{equation}
    \label{prob_comp}
\min_x f(x)=g(x)+h(x)=\mathbb{E}[G(x,\xi)]+\mathbb{E}[H(x,\xi)],
\end{equation}
where $g$ and $h$ are convex.
The composite optimization is a common model for supervised machine learning with regulization,
alternatively nonsmooth convex function or smooth convex function.
And  an algorithm with optimal complexity use the stochastic gradient method for smooth function
 and proximal point method for nonsmooth convex function.

We will analysis the following algorithm.
The first step is Nesterov acceleration
, the second step is a proximal Robbins-Monro gradient for function $g$
, and the third step is stochastic gradient step of function $h$, namely prox-RM-ssgd method .
The second and third steps could be seen as a kind of alterative direction method.

\textsc{\begin{algorithm}[!htbp]
    \caption{prox-RM-ssgd} \label{algo2}
    \begin{algorithmic}[1]{
    \Require Step size $\{\alpha_k\}$, momentum size $\{\theta_k\}$,
    initial value $v_1$, $v_2$, and iteration number $N$,
    \For {$n=1,2,\cdots$}
    Calculate
    $$x_{k+1}=(1+\theta_k)v_{k+1}-\theta_k v_k$$
    and the alterative steps
    $$v_{k+\bullet } \in x_k -a_k \partial G(v_{k+\bullet},\xi)$$
    $$v_{k+1} \in v_{k+\bullet } -a_k \partial H(v_{k+\bullet},\xi)$$
    \EndFor
    }
    \end{algorithmic}
    \end{algorithm}}

    Combined with Lemma \ref{lemma_couple} and Lemma \ref{Lemma_prox1},
    we have the following extension version for composite optimization \ref{prob_comp},
    \begin{lemma}
        \label{lemma_comp}
        Consider the  positive stochastic sequence $r_k$, $\eta_k$, $\zeta_k$, $\rho_k$.
        And sequence $\theta_k\in (0,1)$ is bounded. $a_k$ is a positive decreasing sequence.
        $\sum_{k=1}^\infty \beta_k<+\infty$ , $\sum_{k=1}^\infty \bar{\beta}_k<+\infty$, and
        $$\mathbb{E}[r_{k+2}|\mathcal{F}_{k+1}] \leq (1+\theta_k) r_{k+1}-\theta_k r_k- \eta_k +\beta_k+ hp_k z_{k+1}, k\geq 1$$
        $$\mathbb{E}[z_{k+2}|\mathcal{F}_{k+1}] \leq (1-p_k) z_{k+1} - a_k(\rho_k-\rho_{k-1})+\bar{\beta}_k, k\geq 1.$$
        Then $r_k$ converges to some finite random variable $r_\infty$,a.s
         and  $\sum_{k=1}^\infty \eta_k<+\infty$,a.s.
\end{lemma}

\subsection{Almost surely convergence}
\begin{theorem}
    Consider $\{x_{k}\}$ is generated by Algorithm \ref{algo2}.
    The Assumption \ref{asm_1} holds and momentum parameter $\{\theta_k\}$ is nonincreasing.
    Then $\{x_k\}$ converges to some optimal $x^*$ for problem \eqref{prob_comp}.
\end{theorem}
\begin{proof}
According to the scheme of Nesterov acceleration,
$$\|x_{k+1}-x^*\|^2 = (1+\theta_k) \|v_k-x^*\|^2 - \theta_k \|v_{k-1}-x^*\|^2+ \theta_k(1+\theta_k) \|v_k-v_{k-1}\|^2 .$$
$$\|v_{k+1}-x^*\|^2 = \|v_{k+\bullet}-x^*\|^2-2\alpha_k \langle \tilde{\nabla} H(v_{k+\bullet},\xi_k), v_{k+\bullet}-x^*\rangle +\alpha_k^2M,$$
where
$$\|v_{k+\bullet}-x^*\|^2 = \|x_k-x^*\|^2-2\alpha_k \langle \tilde{\nabla} G(v_{k+\bullet},\xi_k), v_{k+\bullet}-x^*\rangle -\|v_{k+\bullet}-x_k\|^2.$$
Set $r_k=\|v_{k}-x^*\|^2$, both sides take conditional expectation on $\mathcal{F}_k$

 $$\mathbb{E}[r_{k+1}|\mathcal{F}_k]\leq (1+\theta_k)r_k-\theta_kr_{k-1} -2\alpha_k(f(v_{k+\bullet})-f(x^*))+\alpha_k^2M-\|v_{k+\bullet}-x_k\|^2.$$

   $$\|v_{k+1}-v_k\|^2 = \|v_{k+\bullet}-v_k\|^2 -2\alpha_k\langle \tilde{\nabla} F(v_{k+\bullet},\xi_k)\rangle+ \alpha_k^2M-\|v_{k+\bullet}-v_k\|^2$$
   and
   $$\|v_{k+\bullet}-v_k\|^2 = \|x_k-v_k\|^2-2\alpha_k \langle \tilde{\nabla} G(v_{k+\bullet},\xi_k), v_{k+\bullet}-v_k\rangle -\|v_{k+\bullet}-v_k\|^2,$$
Where $\|x_k-v_k\|^2=\theta_k^2\|v_k-v_{k-1}\|^2$.  Set $z_k=\|v_k-v_{k-1}\|^2$, both sides take conditional expectation on $\mathcal{F}_k$.
Then $$\mathbb{E}[z_{k+1}|\mathcal{F}_k]\leq \theta_k^2 z_k -2\alpha_k(f(v_{k+\bullet})-f(v_k))+\alpha_k^2M-\|v_{k+\bullet}-v_k\|^2.$$
According to Lemma \ref{lemma_comp}, $z_k$ converges to 0 almost surely. Similar to Theorem \ref{thm:PRM}, $\{x_k\}$ and $\{v_k\}$ converge to some optimal $x^*$ almost surely.
\end{proof}

\section{Nmerical experiments on Nesterov accelerated methods}

In this section, we consider three problems, the linear least square problem with SGD, the linear least absolute problem with SGD method and prox-RM method, furthermore the Lasso problem with SGD-prox-RM method. There are many famous results for these methods. Here we only give a series of numerical experiments for almost surely convergence.

\subsection{Linear least square problem}
$A\in \mathbb{R}^{m\times n}$, $b\in \mathbb{R}^m$, where $n=20$, $m=2000$, consider the linear least absolute problem $\min\limits_x \sum\limits_{k=1}^m (a^\top_ix-b_i)^2$. Then we take the random index on $\{1,\cdots,m\}$ discrete uniform distribution. Take the random seed 10 of rand in matlab, for example rand('seed',10). Take $v$ with $n\times 1$ uniform distribution on $[0,1]$, and $A$ is a standard normal random matrix multiplied by $I+vv^\top$. $b=Ax_0$, where $x_0$ is the optimal. Step size $\alpha_k=\frac{1}{16(k+3)^{8/9}}$, momentum parameters $\theta$ are constant .

\begin{figure}[!hbtp]
\centering
\label{fig_lsq}
\begin{tabular}{c c}
\includegraphics[width=160pt]{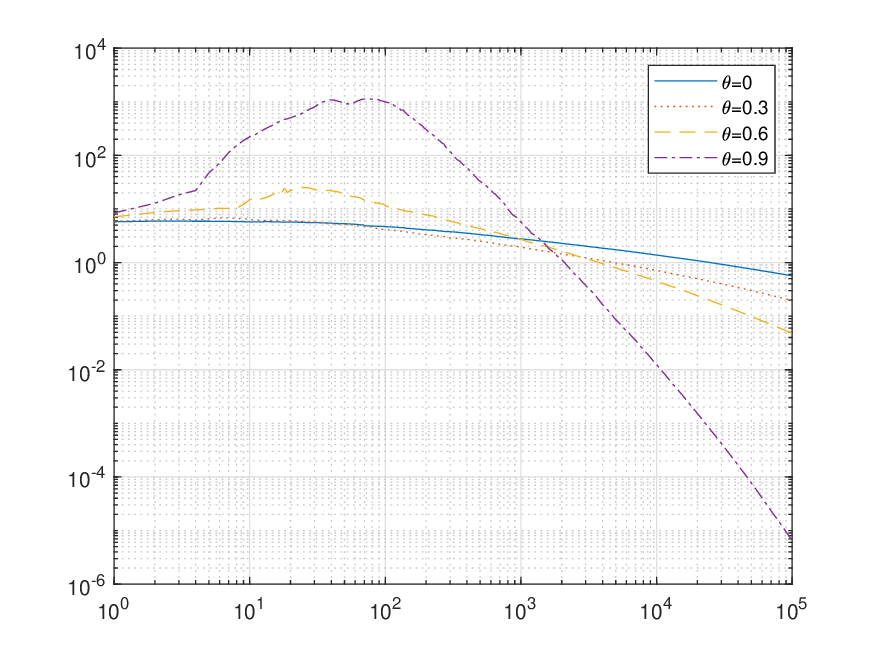}&
\includegraphics[width=160pt]{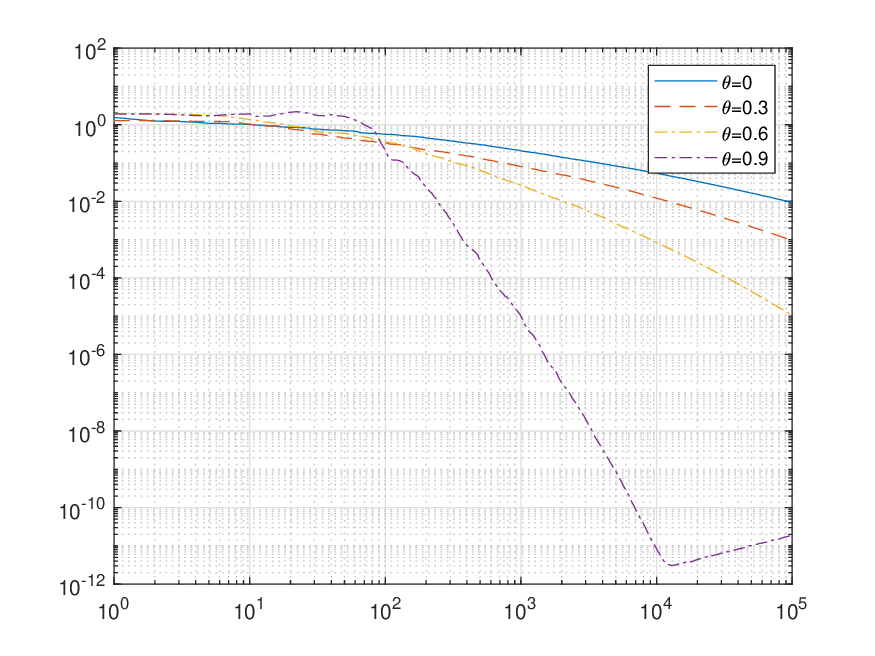}\\
(a) the ssgd method & (b) the prox-RM method
\end{tabular}
\caption{The log-log plot of the stochastic sequence $\{\|v_{k+1}-x^*\|\}$ with the same step size but different momentum parameters for least square problem. }
\end{figure}
The convergence performance is displayed as follows \ref{fig_lsq}. For such a strongly convex function with Lipschitz continuous gradient problem, the nesterov-accelerated ssgd is totally better than the one without, see Figure \ref{fig_lsq} (a). And the prox-RM is more stable than the ssgd method, where there is no gradient exploding process. Also the convergence performance with nesterov-accelerated prox-RM is better than the one without acceleration, see \ref{fig_lsq} (b).

\subsection{Linear least absolute problem}

$A\in \mathbb{R}^{m\times n}$, $b\in \mathbb{R}^m$, where $n=100$, $m=10000$, consider the linear least absolute problem $\min\limits_x \sum\limits_{k=1}^m |a^\top_ix-b_i|$. Then we take the random index on $\{1,\cdots,m\}$ with discrete uniform distribution. Take the random seed 10 of rand in matlab, for example rand('seed',10). Take $v$ with $n\times 1$ uniform distribution on $[0,1]$, and $A$ is a standard normal random matrix multiplied by $I+vv^\top$. $b=Ax_0$, where $x_0$ is the optimal. Step size $\alpha_k=\frac{1}{2(k+3)^{8/9}}$, momentum parameters $\theta$ are constant .

\begin{figure}[!hbtp]
\centering
\label{fig_abs}
\begin{tabular}{c c}
\includegraphics[width=160pt]{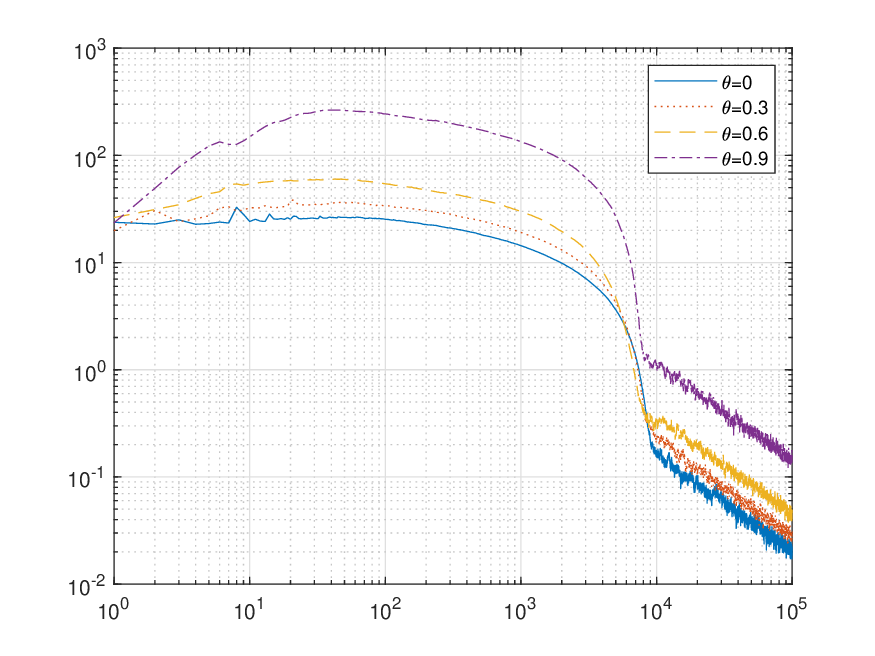} &
\includegraphics[width=160pt]{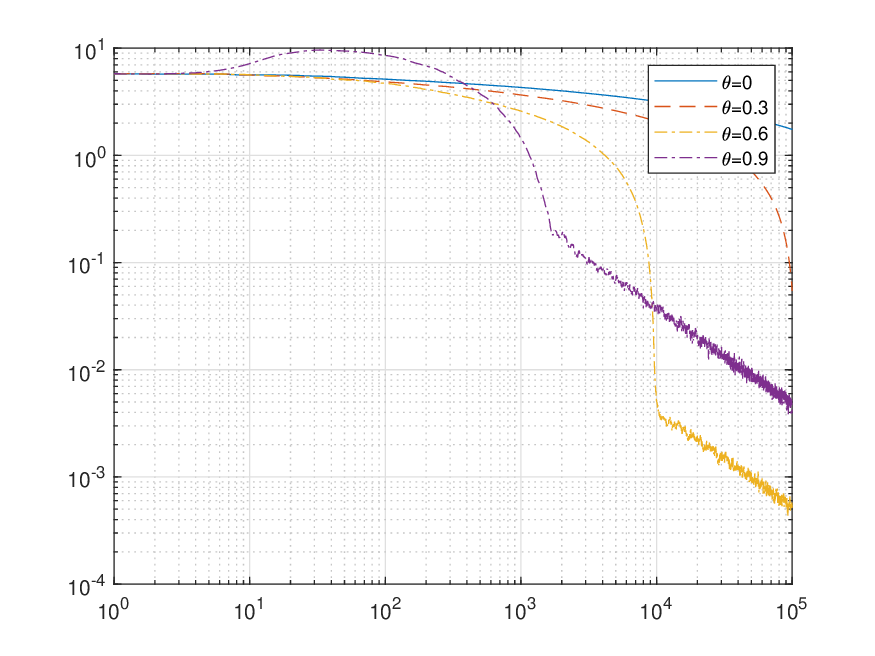}\\
(a) ssgd method & (b) prox-RM method
\end{tabular}
\caption{The log-log plot of the stochastic sequence $\{\|v_{k+1}-x^*\|\}$, from stochastic subgradient method with the same step sequence but different momentum parameters for least absolute problem.}
\end{figure}
The convergence performance is displayed in Figure \ref{fig_abs}.
Although the Nesterov Acceleration method is not better than without the momentum in nonsmooth problem, the almost surely convergence still holds. Also, the prox-RM method is more stable than the ssgd method, where there is no gradient exploding process. Then consider the proximal Robbins-Monro method. Step size $\alpha_k=\frac{1}{4(k+3)^{8/9}}$, momentum parameters $\theta$ are constant.

\subsection{Composite optimization}

Consider the Lasso problem:
\begin{equation}
\min_x \frac{1}{n}\sum_{i=1}^n (a_i^{\top} x-b_i)^2 +\|x\|_1,
\end{equation}
where $A\in \mathbb{R}^{m\times n}$, $b\in \mathbb{R}^m$, $m=10000$, $n=100$. Take the random seed 10 of rand in matlab, for example rand('seed',10). $A$ is a standard normal random matrix and $b$ is a standard normal random vector. Here we only choose the index from $\{1,\cdots,N\}$ randomly and consider the $\|x\|_1$ as a certain function.  The step-size of gradient and proximal is the same $\alpha_k=\frac{1}{20(k+3)^{8/9}}$.

\begin{figure}[!hbtp]
\centering
\label{fig_lasso_ssg}
\includegraphics[width=200pt]{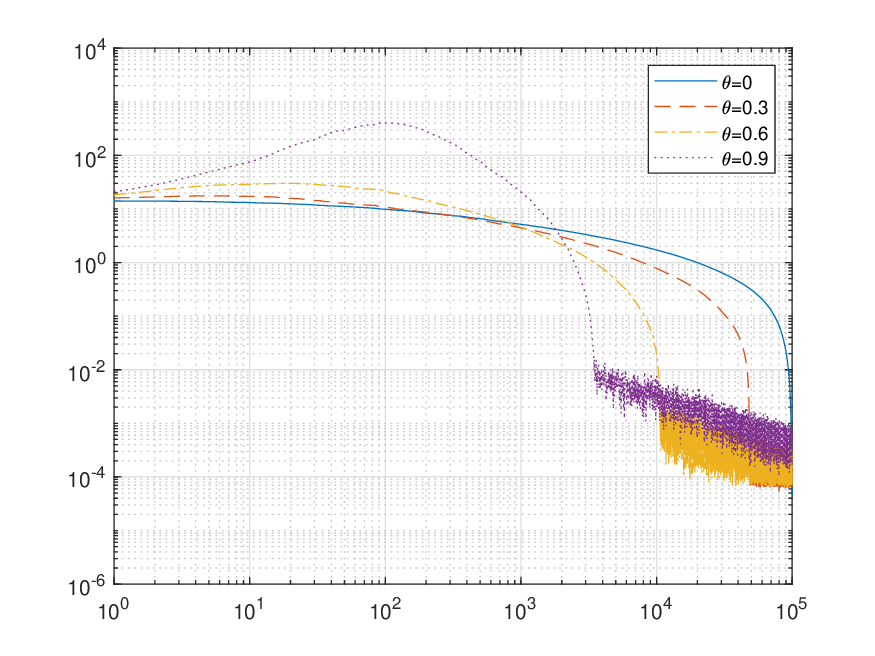}
\caption{The log-log plot of the stochastic sequence $\{\|v_{k+1}-x_N\|\}$, $N=1e5$, different momentum parameters for least absolute problem.}
\end{figure}

The convergence performance is in Figure \ref{fig_lasso_ssg}.

\newpage

\section{Conclusion}

This paper introduces a novel framework for analyzing the convergence of stochastic optimization algorithms, particularly those employing Nesterov accelerated methods.  The key contributions of the paper are twofold:

1. Supermartingale with delayed information: The paper extends the analysis of stochastic sequences to include delayed term, which is a more realistic representation of the stochastic nature of many optimization problems. By incorporating delayed noise into the expected inequalities, the framework captures the temporal aspect of the stochastic environment, providing a more robust and accurate understanding of the algorithm's behavior.

2. Nesterov Accelerated Stochastic Approximation: The paper demonstrates the applicability of the framework to the almost sure convergence of Nesterov accelerated stochastic approximation, a powerful optimization technique. This application highlights the practical significance of the theoretical results, as it ensures that the algorithms will converge to a solution with probability one for both stochastic subgradient method and proximal Robbins-Monro method.

In conclusion, the paper offers a novel and comprehensive framework for analyzing the almost sure convergence of Nesterov accelerated methods.  These findings contribute to the development of more efficient and reliable optimization techniques, particularly in the context of machine learning, data analysis, and control systems.

\bibliographystyle{plain}
\bibliography{ref}

\end{document}